\documentclass[11pt]{article}
\usepackage{amsmath,amsthm,amssymb}
\newtheorem{thm}{Theorem}

\newcommand{\eps}{\varepsilon}

\newcommand{\dist}{{\mathrm{dist}}\,}

\newcommand{\R}{{\mathbb{R}}}
\newcommand{\Z}{{\mathbb{Z}}}

\newcommand{\Q}{{\mathbb{Q}}}
\newcommand{\N}{{\mathbb{N}}}
\newcommand{\calL}{{\mathcal{L}}}

\newcommand{\bomega}{\boldsymbol\omega}
\newcommand{\bOmega}{\boldsymbol\Omega}
\newcommand{\be}{\mathbf{e}}
\newcommand{\bj}{\mathbf{j}}
\newcommand{\ba}{\mathbf{a}}
\newcommand{\bx}{\mathbf{x}}
\newcommand{\by}{\mathbf{y}}
\newcommand{\bv}{\mathbf{v}}
\newcommand{\bq}{\mathbf{q}}
\newcommand{\bp}{\mathbf{p}}
\newcommand{\bM}{\mathbf{M}}

\title{Chaotic dynamics for the two body problem on a sphere}
\author{Sergey Bolotin\\
Steklov Mathematical Institute of Russian Academy of Sciences\\ Moscow, Russia}

\begin{document}

\maketitle

\centerline{Dedicated to Serge Aubry at his 80th birthday}

\begin{abstract}
We prove the existence of chaotic trajectories for the two body problem on a sphere.
The trajectories we construct encounter near-collisions and are similar to the second species solutions
of Poincar\'e  of the classical 3 body problem.
The construction  uses a general result on Lagrangian systems with Newtonian singularities of the potential
which is based on the method of anti-integrable limit of Serge Aubry.
\end{abstract}

\section{Introduction}

Celestial mechanics on surfaces of constant curvature (the sphere and the Lobachevsky plane) goes back to the middle of 19th century (see the history in  \cite{AKN,Bor-Mam1,Diacu}),
and it attracted considerable attention in recent years (see, e.g., \cite{Bor-Mam1,Diacu,Koz-Har,Peres,Shep1,Shep3}).
While it shares many similarities with celestial mechanics in Euclidean space, there are significant differences.
For instance, the two-body problem in Euclidean space reduces to the Kepler problem and is therefore integrable.
The Kepler problem on surfaces of constant curvature closely resembles its Euclidean counterpart: its bounded trajectories are periodic (Serret, 1890).
However, for nonzero curvature, the two-body problem cannot be reduced to the Kepler problem. Numerical evidence (see, e.g., \cite{Bor-Mam1}) suggests that the two-body problem on a sphere exhibits chaotic behavior in certain regions of phase space, although, to the best of our knowledge, this has never been proven analytically.

Several results establish nonintegrability, but only in the meromorphic sense (see e.g.\ \cite{Maci,Shep2,Ziglin}).
However, meromorphic nonintegrability is a much weaker notion than real-analytic nonintegrability and has limited relevance to the system's real dynamics. Indeed, most linear time periodic systems with meromorphic time dependence are meromorphically non-integrable, yet they are completely integrable in the real sense, possessing a full set of quadratic first integrals.

In this paper, we demonstrate chaotic behavior and, in particular,
prove real-analytic nonintegrability of the two-body problem on the sphere in the case where one body's mass
is much smaller than the other's.
Then the heavier body moves nearly inertially, and, to first approximation, the motion of the lighter body is described by the restricted two-body problem.
The proofs in this paper rely on general results in \cite{Bol-Mac:plane}
concerning Lagrangian systems with Newtonian singularities.

The chaotic trajectories we construct pass close to collisions between the bodies. These trajectories resemble the so-called "second species solutions,'' whose existence Poincar\'e claimed to prove in \cite{Poi} for the planar three body problem with masses
$\eps m_1,\eps m_2,1$, where $0<\eps\ll 1$. Second species solutions are those which, in the limit $\eps\to 0$, converge to chains of collision orbits -- what Poincar\'e referred to as "orbites \'a chocs.'' However, his sketch proof of their existence is incomplete, and, as far as we know, rigorous proofs of Poincar\'e's claims have never been published.

For the planar circular restricted three body problem, numerous results on second species solutions have been obtained by astronomers; see the survey in \cite{Henon}. Mathematical contributions include works by Alexeyev \cite{Alex}, Perko \cite{Per}, Henrard \cite{Henrard}, Bruno \cite{Bruno}, Gomez and Olle \cite{Oll-Gom}, and Marco and Niederman \cite{Mar-Nid}. Chaotic second species solutions were constructed in \cite{Bol-Mac:plane} for the planar circular restricted problem and in \cite{Bol-Mac:spat} for the spatial restricted three body problem. In \cite{Mac:magnetic}, MacKay and Pinheiro constructed second species solutions describing the motion of charged particles in a magnetic field.

\section{Main result}

Without loss of generality let the radius of the sphere be one,
then $S^2=\{\bq\in\R^3:|\bq|=1\}$.
The Lagrangian of the 2-body problem on the sphere $S^2$  is
$$
L=\frac12 m_1|\dot \bq_1|^2+\frac 12 m_2|\dot \bq_2|^2+U(\bq_1,\bq_2),
$$
where $U$ is the interaction potential. The Hamiltonian on $T^*(S^2\times S^2)$ is
$$
H=\frac1{2 m_1}|\bp_1|^2+\frac 1{2 m_2}|\bp_2|^2-U(\bq_1,\bq_2),
$$
with the standard symplectic structure on
$$
T^*S^2=\{(\bq,\bp):|\bq|=1,\langle \bp,\bq\rangle =0\}.
$$

The natural analog of the Newtonian gravitational potential for the  sphere is
$$
U(\bq_1,\bq_2)=\frac{k}{\tan\rho},\qquad \rho=d(\bq_1,\bq_2),
$$
where the distance is with respect to the standard metric on $S^2$:
$$
\rho=\arccos\langle \bq_1,\bq_2\rangle.
$$
Explicitly
\begin{equation}
\label{eq:U}
U(\bq_1,\bq_2)=\frac{k \langle \bq_1,\bq_2\rangle}{\sqrt{1-\langle \bq_1,\bq_2\rangle^2}}.
\end{equation}
The potential is singular on $\Delta=\{(\bq_1,\bq_2)\in S^2\times S^2:\bq_1=\pm \bq_2\}$.
As $\rho\to 0$, we have
\begin{equation}
\label{eq:ass}
U\sim \frac{k}{\rho},
\end{equation}
as for the Newtonian potential. Usually one sets $k=\gamma m_1 m_2$.

There are two reasons why the potential (\ref{eq:U}) is the natural analog of the Newtonian potential:

\begin{itemize}
\item
For fixed $\ba\in S^3$, the function $V(\bq)=U(\bq,\ba)$ on the 3-dimensional sphere
satisfies the Laplace-Beltrami equation, just as the Newtonian potential on $\R^3$ (Schering 1873).
\item
For the Kepler problem on a sphere (when one of the particles is fixed),
the motion of the second particle in the potential field $V(\bq)$ is periodic (Serret 1890),
exactly as for the Kepler problem in the Euclidean space.
Darboux (1884) proved a generalization of Bertrand's theorem:
the only other such potential is the analog $k\tan^2\rho$ of the Hooke potential.
\end{itemize}

However, for the results below it is only important that $U$ is $SO(3)$ invariant and has asymptotics (\ref{eq:ass})
as $\rho\to 0$.

Contrary to the 2-body problem in $\R^2$, the  2-body problem on $S^2$ is not integrable: the integrals
of energy $H$ and the angular momentum
$$
\bM=\bq_1\times \bp_1+\bq_2\times \bp_2
$$
corresponding to the $SO(3)$ symmetry are insufficient for integrability.
For fixed angular momentum $\bM=\bM_0\ne 0$,  only the $SO(2)$ symmetry of rotation around $\bM_0$ remains.
Hence the reduced Hamiltonian system has two degrees of freedom.
Explicit reduction was first performed by Shchepetilov \cite{Shep1,Shep2},
and then a simpler reduction was done in  \cite{Bor-Mam1,Bor-Mam3}.

Numerical evidence (see e.g.\ \cite{Bor-Mam1}) shows that the 2-body problem on $S^2$
is chaotic in certain regions of the phase space
but, as far as we know, this was never proved analytically.
There are several results establishing nonintegrability, but only in the meromorphic sense \cite{Maci,Shep2,Ziglin}
which has limited relation  to real dynamics of the system.

We establish chaotic behavior and, in particular, real-analytic nonintegrability
of the two body problem on $S^2$ in the case when one mass is much smaller than the other.

\begin{thm}
\label{thm:2body}
Fix constants $-1/2<a<b$ and
let $\eps_0>0$ and $\sigma_0>0$ be sufficiently small.
If $m_2/m_1<\sigma_0$,
  the 2-body problem on $S^2$ has a chaotic invariant set
on the level sets $\{H=h,\bM=\bM_0\}$ of the first integrals for
$$
a<\frac{m_1^2h}{m_2M_0^2}<b,\quad \frac{m_1^2k}{m_2 M_0^2}<\eps_0.
$$
\end{thm}

In many respects celestial mechanics on a sphere is similar to celestial mechanics on a Lobachevsky space,
see e.g.\ \cite{Bor-Mam3,Diacu}.
In particular, meromorphic nonintegrability still holds, see
\cite{Maci,Shep2}. But an analog of Theorem \ref{thm:2body}
will not hold: for similar values of the first integrals trajectories will escape
to infinity.

As mentioned in the introduction, chaotic trajectories whose existence is claimed in Theorem \ref{thm:2body}
encounter near collisions and are similar in nature to the
``second species solutions''  of Poincar\'e   \cite{Poi} of the Euclidean   3 body problem.

\section{Restricted 2 body problem}

The proof of Theorem \ref{thm:2body} is based on the reduction to the restricted 2 body problem on $S^2$.
We recall a formal derivation. In the formal limit  $m_1\to +\infty$,
the light particle does not influence the motion of the heavy one:
the heavy particle performs the inertial motion on $S^2$.
Hence $\bq_1$ will move with constant angular velocity $\bOmega=\Omega \be$, $|\be|=1$, along a geodesic (great circle)
\begin{equation}
\label{eq:equator}
\Gamma=\{\bq\in S^2:\langle \bq,\be\rangle=0\}.
\end{equation}
Let $\bq_1(0)=\ba\in\Gamma$, then
$$
\bq_1(t)=R_{t\bOmega}\ba=\ba\cos\Omega t+\bj\sin\Omega t,\qquad \bj=\be\times\ba,
$$
where
$R_{t\bOmega}\in SO(3)$
is the  rotation around $\be$ by the angle $t\Omega$.

The motion of the light particle $\bq_2$ is then described by the
time-periodic Lagrangian
$$
L_2(\bq_2,\dot \bq_2,t)=\frac 12m_2|\dot \bq_2|^2+U(\bq_1(t),\bq_2).
$$
If we pass to a rotating coordinate frame, i.e.\ make the change of variables
$\bq_2=R_{t\bOmega}\bx$ and divide the Lagrangian by $m_2$,
we obtain an autonomous  Lagrangian system with the configuration space $S^2$ and the Lagrangian
$$
L(\bx,\dot \bx)=\frac 12|\dot \bx+\bOmega\times \bx|^2+\mu V(\bx),\quad V(\bx)=\frac{\langle \ba,\bx\rangle}{\sqrt{1-\langle \ba,\bx\rangle^2}},
$$
where $\mu=k/m_2$. The potential is singular for  $\bx=\pm\ba$.
Of course this reduction to the restricted problem
is not rigorous; later on it will be made precise.

The restricted 2 body problem  has the Jacobi integral
$$
E(\bx,\dot \bx)=\frac12|\dot \bx|^2-\frac12|\bOmega\times\bx|^2-\mu V(\bx).
$$
The corresponding Hamiltonian on $T^*S^2$ is
$$
H(\bx,\bp)=\frac 12|\bp|^2 -\langle \bp,\bOmega\times \bx\rangle -\mu V(\bx),\quad  \bp=\dot \bx+\bOmega\times \bx.
$$

For $\Omega=0$ we obtain the integrable Kepler problem on $S^2$, see the history in \cite{AKN,Bor-Mam1}.
For $\Omega>0$, the restricted 2 body problem on $S^2$ was studied in several papers,
see e.g.\ \cite{Bor-Mam2,Ziglin,Shep2,Maci}.
Numerical evidence \cite{Bor-Mam1} shows that it exhibits chaotic behavior,
but, as far as we know,
this was never proved analytically.
Non-existence of additional meromorphic first integrals was proved in \cite{Ziglin,Maci}.

We will prove  chaotic behavior of the restricted 2-body problem on $S^2$
for small $\mu/\Omega^2$ and bounded Jacobi constant.

\begin{thm}
\label{thm:main}
Let $-1/2<a<b$, and let $\eps_0>0$ be sufficiently small.
Then for $\mu/\Omega^2<\eps_0$ the restricted 2 body problem on $S^2$
has  a chaotic invariant set on the level set $\{E=h\}$
of the Jacobi integral for $a\Omega^2<h<b\Omega^2$.
\end{thm}

More precisely, an appropriate Poincar\'e map has a hyperbolic invariant set on which it is conjugate
to a topological Markov chain (a subshift of finite type).  In particular, this implies analytic non-integrability \cite{AKN}.
A description of chaotic trajectories is given in Theorem \ref{thm:final}.

In the next section we recall the result in \cite{Bol-Mac:plane} which is the basis for the proof
of Theorem \ref{thm:main}. A more detailed version of Theorem \ref{thm:main} is proved in section \ref{sec:proof}.
In the last section we show how Theorem \ref{thm:main} implies  Theorem \ref{thm:2body}.

\section{Lagrangian systems with singularities}
\label{sec:general}

Let $Q$ be a surface\footnote{In \cite{Bol-Mac:plane} the result is proved also for $\dim Q=3$,
but we need it for $\dim Q=2$ only.} and
\begin{equation}
\label{eq:L}
L(q,\dot q)=\frac 12\|\dot q\|^2+\langle w(q),\dot q\rangle+W(q),
\end{equation}
a smooth\footnote{$C^4$ is certainly enough.} Lagrangian on $TQ$.
Here $\|\;\|^2=\langle\;,\;\rangle$ is a Riemannian metric on $Q$.
Let $V$ be a smooth function on $Q\setminus\Delta$ having
Newtonian singularities on a finite set $\Delta\subset Q$.
This means that in a neighborhood of any point $a\in\Delta$,
\begin{equation}
\label{eq:sing}
V(q)=\frac {f(q)}{\dist(q,a)},\qquad f(a)\ne0,
\end{equation}
where $f$ is a smooth function in a neighborhood of $a$.
The distance is defined by the Riemannian metric.
If $f(a)>0$, the singularity is attracting, and if $f(a)<0$ -- repelling.

Consider a perturbed Lagrangian system
with the configuration space $Q\setminus \Delta$ and the Lagrangian
\begin{equation}
\label{eq:Leps}
L_\eps(q,\dot q)=L(q,\dot q)+\eps V(q).
\end{equation}
Let
$$
E_\eps=E-\eps V,\qquad E(q,\dot q)=\frac 12\|\dot q\|^2-W(q),
$$
be the energy (Jacobi) integral.
We fix $h$ such that the domain $D=\{q \in Q : h+W(q)>0\}$ contains
$\Delta$ and study the system  on the energy level
$\{E_\eps=h\}$.

We call a trajectory $\gamma:[0,\tau]\to Q$
of the unperturbed system a {\it collision orbit} if
$a=\gamma(0),b=\gamma(\tau)\in\Delta$ and $\gamma(t)\notin\Delta$ for $0<t<\tau$ (no early collisions).
We call a collision orbit $\gamma$ with energy $E=h$   {\it nondegenerate} if it is
a nondegenerate critical point of the Maupertuis-Jacobi action functional
$$
J(x)=\int_0^\tau  (2\sqrt{h+W(x)}\|\dot x\|+\langle w(x),\dot
x\rangle)\,dt
$$
on the set of curves $x:[0,\tau]\to D$ joining $a$ with $b$.
The functional  is independent of the parametrisation,
so nondegeneracy is  modulo time reparametrization.
In other words the end points of $\gamma$
are non-conjugate for fixed energy $E=h$.
Then any points $x\in B_\delta(a)$ and $y\in B_\delta(b)$ in small balls around the end points are joined by a unique
trajectory $\gamma_{x,y}:[0,\tau(x,y)]\to Q$ with $E=h$ smoothly depending on $(x,y)$. The action
\begin{equation}
\label{eq:calL}
\calL(x,y)=J(\gamma_{x,y})
\end{equation}
is a smooth function on $B_\delta(a)\times B_\delta(b)$.

To verify the nondegeneracy condition, the following reformulation is useful.
 Let $q(t)=g(a,v,t)$ be the solution of the unperturbed system
with the initial condition $q(0)=a$, $\dot q(0)=v\in T_aQ$.
Collision orbits with $E=h$ connecting $a$ and $b$
correspond to solutions of the equations
\begin{equation}
\label{eq:g}
g(a,v,\tau)=b, \quad E(a,v)=h
\end{equation}
for the variables $v\in T_aQ$ and $\tau>0$.
Then $\gamma$ is a nondegenerate collision orbit
if the pair $v=\dot\gamma(0)$, $\tau$ is a
solution of equations (\ref{eq:g}) at which the rank of the Jacobi matrix is maximal.
Equivalently, the equations
$$
\delta g=\partial_vg(a,v,\tau)\,\delta v+\partial_\tau g(a,v,\tau)\,\delta \tau=0,\quad
\delta E=\partial_v E(a,v)\,\delta v=0
$$
imply $\delta v=0$, $\delta\tau=0$. Here
$$
\partial_\tau g(a,v,\tau)=\dot\gamma(\tau),
\quad
\partial_v E(a,v)\,\delta v=\langle v,\delta v\rangle.
$$

Suppose there exists a finite collection of nondegenerate collision orbits
$\gamma_k:[0,\tau_k]\to D$, $k\in K$, with energy $h$.
We assume that if $\gamma_k(\tau_k)=\gamma_l(0)$, then  $\dot\gamma_k(\tau_k)\ne\pm\dot\gamma_l(0)$
(the changing direction condition).
Let  $\Gamma$ be an oriented  graph with vertices $\Delta$ and edges $\{\gamma_k\}_{k\in K}$.
A bi-infinite path in $\Gamma$ is a  chain of collision orbits $(\gamma_{k_i})_{i\in\Z}$ such that
$\gamma_{k_i}(\tau_{k_i})=\gamma_{k_{i+1}}(0)$.

\begin{thm}
\label{thm:chain}
There exists  $\eps_0>0$ such that
for all $\eps\in(0,\eps_0]$ and any collision chain $(\gamma_{k_i})_{i\in\Z}$:
\begin{itemize}
\item
There exists a unique (up to a time shift) trajectory $\gamma:\R\to D\setminus \Delta$
of energy $h$   shadowing the chain $(\gamma_{k_i})_{i\in\Z}$.
\item
There exists an increasing sequence $(t_i)_{i\in\Z}$ such that
\begin{equation}
\label{eq:upper}
d(\gamma(t),\gamma_{k_{i}}([0,\tau_{k_{i}}]))\le C\eps \quad\mbox{for}\quad t_i\le t\le
t_{i+1}.
\end{equation}
\item
$\gamma$ avoids collisions by a distance of order $\eps$:
\begin{equation}
\label{eq:lower}
\dist(\gamma(t),\Delta)\ge c\eps\quad \mbox{for all}\; t.
\end{equation}
\end{itemize}
\end{thm}

The constants $0<c<C$ depend only on the set $\{\gamma_k\}_{k\in K}$ of collision orbits.
The proof in \cite{Bol-Mac:plane} is based on the method of anti-integrable limit \cite{Aubry}:
trajectories shadowing the collision chain correspond to critical points of the discrete action functional
\begin{equation}
\label{eq:action}
A=\sum_{i\in\Z} \calL_{i}(x_i,y_i)+\eps f_i(y_{i-1},x_i),\qquad x_i\in \partial B_\delta(a_{i-1}),\quad y_i\in \partial B_\delta(a_i),
\end{equation}
where $\calL_i$ is the discrete Lagrangian (\ref{eq:calL}) corresponding to the collision orbit $\gamma_{k_i}$
and $f_i$ a smooth function on a subset in $\partial B_\delta(a_{i-1})\times\partial B_\delta(a_i)$,
$a_i=\gamma_{k_i}(\tau_{k_i})$.
For $\eps\to 0$ the functional splits into a sum of independent functions $\calL_i$,
and each of them has a nondegenerate critical point. Hence we are in a situation of the anti-integrable limit.

Since trajectories in Theorem~\ref{thm:chain} were obtained by the method of anti-integrable limit,
they form a hyperbolic invariant set in  $\{E_\eps=h\}$
on which the system is a
suspension of a subshift of finite type (topological Markov chain).
The Lyapunov exponent is of order $|\ln\eps|$.
The topological entropy is positive if e.g.\ there are at least 2 collision orbits
colliding with the same point in $\Delta$.

\section{Proof of Theorem \ref{thm:main}}
\label{sec:proof}

To get rid of an extra parameter $\Omega$, we introduce new time $\hat t=\Omega t$ and replace the Lagrangian $L$
by $\hat L=L/\Omega^2$. Jacobi constant is replaced by $\hat h=h/\Omega^2$.
Dropping the hats for simplicity, we obtain the Lagrangian
\begin{equation}
\label{eq:tildeL}
L_\eps(\bx,\dot \bx)=\frac 12|\dot \bx+\be\times \bx|^2+\eps V(\bx),\qquad \eps=\mu/\Omega^2,
\end{equation}
Now there are just two parameters: $\eps$ and the Jacobi constant
$$
h=\frac 12|\dot \bx|^2-\frac 12|\be\times \bx|^2-\eps V(\bx).
$$
We fix $h>-1/2$. Let
\begin{eqnarray*}
&A(h)=\Q\cap (1-\sqrt{2h+1},1+\sqrt{2h+1}),\qquad -1/2<h<0,\\
&A(h)=\Q\cap (-1+\sqrt{2h+1},1+\sqrt{2h+1}),\qquad h>0.
\end{eqnarray*}
Take a finite set $K\subset A(h)$. We represent each $\omega\in K$ by an irreducible fraction $\omega=k/n$, $n>0$.

Theorem \ref{thm:main} will be deduced from the following result.

\begin{thm}
\label{thm:final}
Let $\eps_0>0$ be sufficiently small.
For $0<\eps<\eps_0$ and any sequence $(\omega_i=k_i/n_i)_{i\in \Z}\in K^\Z$ there exists a trajectory
$\gamma:\R\to S^2\setminus\Delta$, $\Delta=\{\pm\ba\}$,
with the Jacobi integral $h$ and a sequence $(t_i)_{i\in\Z}$ such that:
\begin{itemize}
\item
$|t_{i+1}-t_i-n_i|\le C\eps$,
and $0<c\eps\le d(\gamma(t_i),\ba_i)\le C\eps$, where the sequence $\ba_i\in\Delta$
satisfies $\ba_{i+1}=(-1)^{k_i+n_i}\ba_i$.
\item
The path $R_{-t\bOmega}\gamma(t)$, $t_i\le t\le t_{i+1}$,
is $O(\eps)$-close to the great circle through $\pm\ba$
with inclination
$$
\alpha_i=\arccos \frac{\omega_i^2/2-h-1}{\sqrt{2h+1}},
$$
to the equator (\ref{eq:equator}) and it makes $\sim[k_i/2]$
revolutions around the great circle between near collisions with $\Delta$.
\item
These trajectories form a hyperbolic invariant set with Lyapunov exponents of order $|\ln\eps|$.
\end{itemize}
\end{thm}

The constant $\eps_0$ depends on $h$ and on the set $K\subset A(h)$. However, for fixed constants
$-1/2<a<b$ and $N>0$, if we take $h\in (a,b)$ and $n\le N$ for every fraction $\omega=k/n\in K$,
then $\eps_0=\eps_0(a,b,N)$ will be independent of $h$ and $K$.
This implies Theorem \ref{thm:main}. This remark will be also needed
for the proof of Theorem \ref{thm:2body}.

\medskip

\noindent{\it Proof.}
We need to check that conditions of Theorem \ref{thm:chain} are satisfied.
Set $\eps=0$. Then the potential $V$ disappears, and the Lagrangian (\ref{eq:tildeL})
describes the geodesic flow on $S^2$ in the rotating coordinate frame.
Hence the trajectories are doubly periodic.

A trajectory starting at $\ba\in \Gamma$ with relative velocity
$\dot \bx(0)=\bv\perp \ba$ is given by
$$
\bx(t)=R_{-t\be}R_{t\bomega}\ba,\qquad \bomega=\ba\times (\bv+\be\times \bv)=\ba\times \bv+\be.
$$
We will use $\bomega$ for the parameter, then
$$
\bv=(\bomega-\be)\times \ba.
$$
The Jacobi constant of the trajectory $\bx(t)$ is
\begin{equation}
\label{eq:h}
h=\frac 12|\bv|^2-\frac 12|\be\times \ba|^2
=\frac 12|(\bomega-\be)\times \ba|^2-\frac 12=\frac12|\bomega|^2-\langle\bomega,\be\rangle.
\end{equation}

Suppose $\langle \bv,\be\rangle\ne 0$, then  the path $R_{t\bomega}\ba$ crosses the equator
$\Gamma=\{\bx\in S^2:\langle \be,\bx\rangle =0\}$ for $t=\tau=\pi k/\omega$, $k\in\N$,
at the point $\ba$ or $-\ba$ depending if $k$ is even or odd: $R_{\tau\bomega}\ba=(-1)^k\ba$.
Then $\bx(\tau)=R_{-\tau\be}R_{\tau\bomega}\ba\in\Delta$ if $t=\pi n$ for $n\in\N$.
More precisely, $\bx(\tau)=(-1)^{k+n}\ba$.

Let $\gamma:[0,\tau]\to S^2$ be a collision orbit starting and ending at a singularity:
$\gamma(0),\gamma(\tau)\in\Delta=\{\pm\ba\}$.
For definiteness suppose $\gamma(0)=\ba$.
Then $\tau=\pi n=\pi k/\omega$, where $\omega=k/n$.
If the fraction is irreducible, there will be no early collisions:
$\gamma(t)\notin\Delta$ for $0<t<\tau$ and $\gamma(\tau)=(-1)^{k+n}\gamma(0)$.

We need the collision orbit to have the given value of the Jacobi constant $h>-1/2$.
Let $\theta$ be the angle between $\bomega$ and $\be$, then
the Jacobi constant equals $h$ when $\cos\theta=\omega/2-h/\omega$.
This equation has a solution $\theta\in (0,\pi)$ for $\omega\in A(h)$.
For this value of $\theta$, the angle $\alpha$ between $\bv$ and $\bj=\be\times\ba$ satisfies
$$
\cos\alpha=\frac{\omega\cos\theta-1}{\sqrt{2h+1}}=\frac{\omega^2/2-h-1}{\sqrt{2h+1}}.
$$
Since this is monotone in $\omega$, for different $\omega$ the changing direction condition will hold.

It remains to check the nondegeneracy condition. For definiteness let $k+n$ be even,
then $\gamma(\tau)=\gamma(0)=\ba$.
Equations (\ref{eq:g}) for the  variables  $\bomega,\tau$ have the form
$$
R_{\tau\bomega}\ba=R_{\tau\be}\ba,\quad \frac12|\bomega|^2-\langle\bomega,\be\rangle =h,\quad \langle \bomega,\ba\rangle=0,
$$
where
$$
R_{\tau\bomega}\ba=\ba \cos(\tau\omega)+\frac{\bomega\times\ba}{\omega}\sin(\tau\omega),\quad
R_{\tau\be}\ba=\ba\cos\tau+\bj\sin\tau.
$$
Hence for $\tau=\pi n=\pi k/\omega$ we obtain
$$
\delta (R_{\tau\bomega}\ba - R_{\tau\be}\ba)|_{\tau=\pi n}=\bomega\times \ba\,\delta(\tau\omega)+\frac{\bomega\times\ba}{\omega}(-1)^k\tau\,\delta\omega-\bj\,\delta\tau.
$$
Taking the dot product with $\bomega$,
we obtain $\delta\tau=0$, and then $\delta\omega=0$.
Since $\langle\delta\bomega,\be\rangle=\langle\delta\bomega,\ba\rangle=0$, we have $\delta\bomega=0$.

Nondegeneracy of collision orbits is proved. Now Theorem \ref{thm:final} follows from Theorem \ref{thm:chain}.
\qed

\section{Proof of Theorem \ref{thm:2body}}
\label{sec:last}

We need to make precise the relation between the 2 body problem and the restricted 2 body problem.
The 2 body problem on $S^2$ has 4 degrees of freedom and $SO(3)$ symmetry with the corresponding angular
momentum integral $\bM$.
For a fixed value of  $\bM=\bM_0=M_0\be$ reducing symmetry yields a Hamiltonian system
on a 4-dimensional symplectic manifold $\{\bM=\bM_0\}/SO(2)$.
This reduction is done in several papers, see e.g.\ \cite{Shep1,Shep2,Bor-Mam1,Bor-Mam2}.
We will use the reduced system obtained in \cite{Shep2}, but  represent it  in a coordinate-free way
as a Hamiltonian system on $T^*S^2$ with the standard symplectic structure and the Hamiltonian
\begin{equation}
\label{eq:Hred}\nonumber
H(\bx,\bp)=\frac{|\bp|^2}{2m}-\frac{\langle \bp,\ba\times \bx\rangle^2}{m_1}
+\frac{\sqrt{M_0^2-\langle\bp,\ba\times\bx\rangle^2}}{m_1}\langle\bp,\be\times\bx\rangle+k V(\bx),
\end{equation}
where $m=m_1m_2/(m_1+m_2)$ is the reduced mass. Set $M_0=m_1\Omega$. Then
\begin{eqnarray*}
H(\bx,\bp)&=&\frac{|\bp|^2}{2m_2} +\Omega\langle\bp,\be\times\bx\rangle+k V(\bx)\\
&& +\frac{1}{m_1}\left(\frac{|\bp|^2}{2}-\langle \bp,\ba\times \bx\rangle^2
-\frac{\langle\bp,\ba\times\bx\rangle^2\langle \bp,\be\times\bx\rangle}
{\sqrt{\Omega^2 m_1^2- \langle\bp,\ba\times\bx\rangle^2}+\Omega m_1}\right).
\end{eqnarray*}
We introduce new momentum $\by=\bp/(m_2\Omega)$ and, scaling time,
replace the Hamiltonian by $\hat H=H/(m_2\Omega^2)$.
We obtain:
\begin{eqnarray}
\nonumber
\hat H(\bx,\by)&=&\frac{|\by|^2}{2} +\langle\by,\be\times\bx\rangle+\eps V(\bx)\\
&&+\sigma\left(\frac{|\by|^2}{2}-\langle \by,\ba\times \bx\rangle^2
-\sigma\frac{\langle\by,\ba\times\bx\rangle^2\langle \by,\be\times\bx\rangle}
{\sqrt{1- \sigma^2\langle\by,\ba\times\bx\rangle^2}+1}\right).
\nonumber
\end{eqnarray}
where $\sigma=m_2/m_1$ and $\eps=k/(m_2\Omega^2)$.
Thus
$$
\hat H(\bx,\by)=H_\eps(\bx,\by)+\sigma F(\bx,\by,\sigma),
$$
where $H_\eps$ is the Hamiltonian of the restricted 2 body problem with the Lagrangian
(\ref{eq:tildeL}), and $F$ is a smooth function on $T^*S^2$ for $|\by|<\sigma^{-1}$.

Let $\Lambda$ be the hyperbolic invariant set on $\{H_\eps=\hat h\}$ constructed
in Theorem \ref{thm:final} for small $\eps>0$. By continuity of compact hyperbolic sets under perturbations,
for small $\sigma_0(\eps)>0$ and $0<\sigma<\sigma_0$ the hyperbolic set $\Lambda$
gives a hyperbolic invariant set of the reduced system on $\{\hat H=\hat h\}$,
and hence a chaotic invariant set of the 2 body problem on the level set $\{H=h,\bM=M_0\be\}$ for
$h=m_2\Omega^2\hat h$ and $M_0=m_1\Omega$. However, since $\sigma_0$ depends on $\eps$,
we do not get Theorem \ref{thm:2body} in this way.

A better estimate for $\sigma$ may be obtained as follows.
Let us assume $-1/2<a<\hat h<b$ and take the set $K$ as in the remark after Theorem \ref{thm:main}.
Then the discrete action functional (\ref{eq:action}) will have uniformly anti-integrable structure.
For the perturbed Hamiltonian $\hat H$ the functional will have the same structure, but now
its terms will smoothly depend on $\sigma$.
This implies that for small but unrelated $\eps_0>0$ and $\sigma_0>0$
the anti-integrability will hold for $0<\eps<\eps_0$ and $0<\sigma<\sigma_0$.
Theorem \ref{thm:2body} will follow.
\qed

\section*{Acknowledgements}

This work was performed at the Steklov International Mathematical Center
and supported by the Ministry of Science and Higher Education of the Russian Federation (agreement no. 075-15-2025-303).

\end{document}